\documentclass[10pt]{svmult}

\usepackage{graphicx}
\usepackage[table,pdftex,hyperref,svgnames]{xcolor}
\usepackage{url,doi}

\usepackage{mathptmx}       % selects Times Roman as basic font
\usepackage{helvet}         % selects Helvetica as sans-serif font
\usepackage{courier}        % selects Courier as typewriter font
\usepackage{type1cm}        % activate if the above 3 fonts are
\usepackage[margin = 1.0in]{geometry}
\usepackage{mathtools}

\usepackage{amsmath}
\usepackage{amsfonts}
\usepackage{subfigure}
\usepackage{bm}
\usepackage{amssymb}
\usepackage{epstopdf}
\usepackage{mathtools}  
\usepackage{tabulary}
\usepackage{booktabs}
\usepackage{cite}

\usepackage{enumitem}

\DeclareSymbolFont{bbold}{U}{bbold}{m}{n}
\DeclareSymbolFontAlphabet{\mathbbold}{bbold}
\newcommand{\vect}[1]{\mathbbold{#1}}

\newcommand{\QED}{\begin{flushright}\qed \end{flushright}}

\newcommand{\intr}{\operatorname{int}}

\def \I{\boldsymbol{I}}

\graphicspath{{./images/}}

\begin{document}

\setlength{\parskip}{0.2cm}
\setlength\parindent{0em}
\setlist{nolistsep}

\chapter*{\Large{LaSalle Invariance Principle for Discrete-time Dynamical Systems:\\ A Concise and Self-contained Tutorial}}

\vspace{-2.5cm}{Wenjun Mei and Francesco Bullo\\
  Center for Control, Dynamical-system, and Computation\\
University of California at Santa Barbara}
\vspace{1.5cm}

\abstract{As a method to establish the Lyapunov stability of
  differential/difference equations systems, the LaSalle invariance
  principle was originally proposed in the 1950's and has become a
  fundamental mathematical tool in the area of dynamical systems and
  control. In both theoretical research and engineering practice,
  discrete-time dynamical systems have been at least as extensively
  studied as continuous-time systems. For example, model predictive
  control is typically studied in discrete-time via Lyapunov
  methods. However, there is a peculiar absence in the standard
  literature of standard treatments of Lyapunov functions and LaSalle
  invariance principle for discrete-time nonlinear systems. Most of
  the textbooks on nonlinear dynamical systems focus only on
  continuous-time systems. For example, the classic textbook by Khalil
  on nonlinear systems~\cite{HKK:02} relegates discrete-time systems
  to a few exercises at the end of Chapter 4. The textbook by
  Vidyasagar~\cite{MV:02} does not present the discrete-time LaSalle
  invariance principle. In Chapter 1 of the book by
  LaSalle~\cite{JPL:76}, the author establishes the LaSalle invariance
  principle for difference equation systems. However, all the useful 
  lemmas in~\cite{JPL:76} are given in the form of
  exercises with no proof provided. In this document, we provide the
  proofs of all the lemmas proposed in~\cite{JPL:76} that are needed
  to derive the main theorem on the LaSalle invariance principle for
  discrete-time dynamical systems. We organize all the materials in a
  self-contained manner. We first introduce some basic concepts and
  definitions in Section 1, such as dynamical systems, invariant sets,
  and limit sets.  In Section 2 we present and prove some useful
  lemmas on the properties of invariant sets and limit sets. Finally,
  we establish the original LaSalle invariance principle for
  discrete-time dynamical systems and a simple extension in
  Section~3. In Section 4, we provide some references on extensions of
  LaSalle invariance principles for further reading. This document is
  intended for educational and tutorial purposes and contains lemmas
  that might be useful as a reference for researchers.}

\vspace{1cm}
\section{Basic Concepts: Dynamical System, Motion, Limit Set, and Invariant Set}

Before reviewing the basic concepts in discrete-time dynamical systems, we first introduce some frequently used notations. Let $\mathbb{N}$ be the set of natural numbers, i.e, $\{0,1,2,\dots\}$. Denote by $\mathbb{Z}$ and $\mathbb{Z}_+$ the set of integers and positive integers respectively. The set of real numbers is denoted by $\mathbb{R}$ and the $m$-dimension Euclidean space is denoted by $\mathbb{R}^m$. Let $\vect{0}_m$ be the all-zeros $m\times 1$ column vector and let $\vect{1}_m$ be the all-ones $m\times 1$ column vector. We use $\phi$ to denote the empty set. For any sequence $\{x_k\}_{k\in \mathbb{N}}$, by $\{x_k\}\to y$ we mean $x_k \to y$ as $k\to \infty$. Without causing any confusion, sometimes we omit the subscript $k\in \mathbb{N}$ when we refer to a sequence $\{x_k\}_{k\in \mathbb{N}}$. Given a set $A$, denote by $|A|$ the cardinality of $A$.

\subsection{Discrete dynamical system}

Given a map $T:\mathbb{R}^m \to \mathbb{R}^m$, the following equations system:
\begin{equation}\label{eq:def-1st-order-diff-eq}
x(n+1)=T\big( x(n) \big), \quad \text{for any }n\in \mathbb{N},
\end{equation}
is referred to as a \emph{$m$-dimension first-order difference equations system}. Equation~\eqref{eq:def-1st-order-diff-eq} together with an additional condition $x(0)=x_0\in \mathbb{R}^m$ defines an \emph{initial-value problem} for the $m$-dimension first-order difference equations system. A sequence $\{x(n)\}_{n\in \mathbb{N}}$ is the solution to the initial-value problem if $x(n)=T^n(x_0)$ for any $n\in \mathbb{N}$. Similarly, for any $k\in \mathbb{Z}_+$, a map $g:\underbrace{\mathbb{R}^m\times \dots \times \mathbb{R}^m }\limits_{k} \to \mathbb{R}^m$ defines a \emph{$m$-dimension $k$-th order difference equations system}, and the corresponding initial-value problem is of the following form:
\begin{equation*}
\begin{cases}
\displaystyle x(n+1) & = g\big( x(n),\dots,x(n-k+1) \big),\\
\displaystyle x(k-1) & = x_{k-1},\\
\displaystyle & \dots  \\
\displaystyle x(0) & = x_0.
\end{cases}
\end{equation*} 

\begin{lemma}[Equivalence to first-order system]\label{lem-transf-first-order}
Any one-dimension $m$-th order difference equations system
\begin{equation*}
u(n+1) = g\big( u(n),\dots,u(n-m+1) \big)
\end{equation*}
is equivalent to an $m$-dimension first-order difference equations system.
\end{lemma}

\textbf{Proof:} For any $n\ge m$, let $x_m(n)=u(n-m+1)$, $x_{m-1}(n)=u(n-m+2)$, $\dots$, $x_2(n)=u(n-1)$, and $x_1(n)=u(n)$. Define $x(n)\in \mathbb{R}^m$ as $x(n)=\big(x_1(n),\dots,x_m(n)  \big)^{\top}$. We thereby have:
\begin{align*}
x_1(n+1) & = g\big( x_1(n),\dots,x_m(n) \big),\\
x_2(n+1) & = x_1(n),\\
         & \dots \\
x_m(n+1) & = x_{m-1}(n).
\end{align*}
Define $T:\mathbb{R}^{m}\to \mathbb{R}^m$ as
\begin{equation*}
T\big( x(n) \big) = 
\begin{bmatrix}
g\big( x_1(n),\dots,x_m(n) \big) \\
x_1(n) \\
\dots \\
x_{m-1}(n)
\end{bmatrix}.
\end{equation*}
We obtain an $m$-dimension first order difference equations system $x(n+1)=T\big( x(n) \big)$. This concludes the proof. \QED

For simplicity, in the rest of this document, whenever we refer to a difference equations system, we assume it is $m$-dimension and first-order, unless specified. In addition, no matter whether it is emphasized, the map $T:\mathbb{R}^m \to \mathbb{R}^m$ is always assumed to be continuous in this document.

\begin{definition}[discrete semi-dynamical system and dynamical system]\label{def:(semi)-dyn-sys}
A map $\pi:\mathbb{N}\times \mathbb{R}^m \to \mathbb{R}^m$ is a discrete semi-dynamical system on $\mathbb{R}^m$ if, for any $n,k\in \mathbb{N}$ and any $x\in \mathbb{R}^m$, 
\begin{enumerate}[label={\arabic*)}]
\item $\pi(0,x)=x$;
\item $\pi\big( n,\pi(k,x) \big) = \pi(n+k,x)$ (semi-group property);
\item $\pi$ is continuous. 
\end{enumerate}
The map $\pi$ is a discrete dynamical system if 2) above holds for any $n,k\in \mathbb{Z}$ (group property).
\end{definition}

\begin{remark}
The semi-group property implies the uniqueness of the solution to equations $x(n+1)=\pi\big(x(n)\big)$, with the initial condition $x(0)=x_0$, in the forward direction, i.e., for any $n\in \mathbb{N}$. The group property leads to the uniqueness of the solution in both direction, i.e., for any $n\in \mathbb{Z}$.
\end{remark}

\subsection{Limit set and invariant set}
\emph{Invariant set} and \emph{limit set} are two important concepts in dynamical systems. Limit set characterizes the asymptotic behavior of a difference equations system, i.e., the limit behavior of $T^n(x_0)$, as $n\to \infty$. A compact invariant set implies the existence of limit sets. In this subsection we present the definitions of invariant set and limit set.

\begin{definition}[Invariant sets]\label{def:inv-set}
Given a map $T:\mathbb{R}^m \to \mathbb{R}^m$ and a set $H\subset \mathbb{R}^m$, define the set $T(H)=\{y\in \mathbb{R}^m\,|\,y=T(x)\text{ for some }x\in H\}$. The set $H$ is positively invariant if $T(H)\subset H$, negatively invariant if $T(H)\supset H$, and invariant if $T(H)=H$. A set $H\subset \mathbb{R}^m$ is invariantly connected if 
\begin{enumerate}[label={\arabic*)}]
\item $H$ is closed and invariant;
\item $H$ is not a union of two non-empty disjoint closed invariant sets.
\end{enumerate}
For any $E\subset \mathbb{R}^m$, a set $M$ is the largest invariant set in $E$, if $M\subset E$, $T(M)=M$, and $T(A\cup M)\neq A\cup M$ for any set $A$ such that $A\subset E$ and $A\cap M = \phi$.
\end{definition}

\begin{remark}\label{remark:inv-set-connected-largest} 
The following statements can be inferred from Definition~\ref{def:inv-set}: 1) Any set of countably many (more than one) isolated fixed points is not invariantly connected; 2) If $M$ is the largest invariant set in $E$, then, for any $A\subset E$, $T(A\cup M)=A\cup M$ implies $A\subset M$.
\end{remark}

Before presenting the definition of limit set, we first introduce some preliminary notions.

\begin{definition}[Motion, periodicity, fixed point, and extension of motion]\label{def:motion}
Given a map $T:\mathbb{R}^m \to \mathbb{R}^m$ and a vector $x\in \mathbb{R}^m$, the motion from $x$ refers to the sequence $\{T^n(x)\}_{n\in \mathbb{N}}$ and is denoted by $T^n x$. A motion $T^n x$ is periodic if there exists some $k\in \mathbb{Z}_+$ such that $T^k(x)=x$. The least such $k$ is referred to as the period of the motion $T^n$. A point $x$ is called a fixed point of $T$ if the motion $T^n x$ has period $k=1$, i.e., $T(x)=x$. A set $T_n x = \{T_n(x)\}_{n\in \mathbb{Z}}\subset \mathbb{R}^m$ is referred to as an extension of the motion $T^n x$ if $T_0 x = x$ and $T(T_n x) = T_{n+1} x$ for any $n\in \mathbb{Z}$.
\end{definition}

\begin{remark}
The extension of a motion $T^n x$ is not unique if the map $T$ is non-invertible.
\end{remark}

\begin{definition}[Distance and convergence]\label{def:dist-conv}
For any point $x\in \mathbb{R}^m$ and set $S\subset \mathbb{R}^m$, define the distance between $x$ and $S$ as $\rho(x,S)=\inf\limits_{y\in S} \lVert x-y \rVert$, where $\lVert \cdot \rVert$ is some norm defined in $\mathbb{R}^m$. A motion $T^n x$ converges to the set $S$ if $\lim\limits_{n\to \infty} \rho\big( T^n(x),S \big)=0$.
\end{definition}

\begin{definition}[Interior, closure, and boundary]\label{def:int-cls-bdr}
For any $x\in \mathbb{R}^m$ and $r>0$, define the open ball around $x$ with radius $r$ as the set $\mathcal{B}(x,r) = \big{\{} y\in \mathbb{R}^m \,\big|\, \lVert y-x \rVert <r \big{\}}$.
For any set $S\subset \mathbb{R}^m$, $x\in S$ is an interior point of $S$ is there exists $\epsilon>0$ such that $\mathcal{B}(x,\epsilon)\subset S$. Denote by $\intr{S}$ the set of all the interior points of $S$. Define the closure of $S$ as the set $\overline{S} = \big{\{} x\in \mathbb{R}^m \,\big|\, \rho(x,S)=0 \big{\}}$. Define the boundary of $S$ as the set $\partial S = \overline{S}\setminus \intr{S}$.
\end{definition}

Now we present the definitions of limit point and limit set.
\begin{definition}[Limit point and limit set]\label{def:lim-set}
Given a motion $T^n x$, $y\in \mathbb{R}^m$ is a limit point of the motion $T^n x$ if there exists a subsequence $\{n_k\}_{k\in \mathbb{N}}$ such that $n_k \to \infty$ and $T^{n_k} (x) \to y$ as $k\to \infty$. When there is no ambiguity about the map $T$, we also refer to $y$ as a limit point of $x$. The set of all the limit points of $x$ is referred to as the limit set of $x$, denoted by $\Omega(x)$. That is,
\begin{equation*}
\Omega(x) = \Big{\{} y\in \mathbb{R}^m \,\Big|\, \text{There exists sequence }\{n_k\}\subset \mathbb{N} \text{ such that }n_k\to \infty\text{ and }T^{n_k}(x)\to y\text{ as }k\to \infty  \Big{\}}.
\end{equation*}
Given a set $H\subset \mathbb{R}^m$, the limit set of the set $H$ is denoted by $\Omega(H)$ and defined as
\begin{equation*}
\Omega(H) = \Big{\{} y\in \mathbb{R}^m \,\Big|\, \text{There exists sequence } \{n_k\}\subset \mathbb{N}\text{ and }\{y_k\}\subset H\text{ such that }n_k \to \infty \text{ and }T^{n_k}(y_k)\to y\text{ as }k\to \infty \Big{\}}.
\end{equation*}
\end{definition}

\section{Properties of Invariant Sets and Limit Sets}
In this section we present and prove some important properties of the limit sets and the invariant sets of difference equations systems. These properties will be used in the proof of the LaSalle invariant principle. 

\subsection{Properties of invariant sets}
We first present some lemmas on the properties of invariant sets.

\begin{lemma}[Invariantly connected set and periodic motion]\label{lem:inv-conn-period}
For any given continuous map $T:\mathbb{R}^m \to \mathbb{R}^m$, suppose $H$ is an invariant set with finite elements. The set $H$ is invariantly connected if and only if $H$ is a periodic motion $T^n x$.
\end{lemma}

\textbf{Proof:} (1) ``if'' part: Suppose $H$ is a periodic motion with period $K$. By definition, $H=\{x,T(x),\dots,T^{K-1}(x)\}$ is invariant and closed. It is also straightforward to check that $H=\{x,T(x),\dots,T^{K-1}(x)\}$ cannot be a union of any two non-empty disjoint closed invariant set. Therefore, $H$ is invariantly connected.

(2) ``only if'' part: We prove the ``only if'' part by constructing a periodic orbit that passes through all the points in $H$. For any given $x\in H$, the motion $T^nx=\cup_{k=0}^{\infty}\{T^k(x)\}$ satisfies 
\begin{align*}
T(T^nx)\subseteq T^nx \subseteq H.
\end{align*}
Now we prove by contradiction that $T^nx = H$. Suppose $T^n x\neq H$. Then $H\setminus T^n x$ cannot be an invariant set, otherwise,
\begin{enumerate}
\item either $T(T^n x)=T^n x$ and $H\setminus T^n x$ are both invariant sets, which contradicts the fact that $H$ is invariantly connected;
\item or $T(T^n x)\subset T^n x$, which leads to
   \begin{align*}
   T(H) = T(H\setminus T^n x) \cup T(T^n x) \subset (H\setminus T^n x) \cup T^n x = H
   \end{align*}
   and contradicts the fact that $H$ is an invariant set.
\end{enumerate}
Since $H\setminus T^n x$ is not invariant, there exists $y\in H\setminus T^n x$ such that $T(y)\in T^n x$. As a consequence, 
\begin{align*}
|T(H)| = \Big| T\big( H\setminus (T^n x\cup \{y\}) \big) \cup T(\{y\}) \cup T(T^n x) \Big| \le \big| H\setminus (T^n x\cup \{y\}) \big| + |T^n x| \le |H|-1< |H|, 
\end{align*}
which contradicts the fact that $H$ is invariant. Therefore, the original assumption that $T^n x\neq H$ does not hold. Now we have proved that, for the given $x$, $T^n x = H$. That is, for any $y\in H$, there exists $k\in \mathbb{N}$ such that $y=T^k (x)$.

Since $H$ is finite, we can denote $H$ by $H=\{x,y_1,\dots, y_m\}$. For any $i\in \{1,\dots,m\}$, define
\begin{align*}
k_i = \min \Big{\{} k\in \mathbb{N}\,|\, T^k(x) = y_i \Big{\}}.
\end{align*}
By definition, $k_i\neq k_j$ if $i\neq j$. Without loss of generality, assume that $k_1<k_2<\dots<k_m$. (This can always be achieved by relabelling $y_1,\dots, y_m$.) Since $x=T^0(x)$, it must hold that $k_1=1$, otherwise $T(x)=x$, which leads to $T^n x = \{x\}\neq H$. Since $x=T^0(x)$ and $y_1 = T^1(x)$, we in turn have $k_2 = 2$, otherwise either $T^2(x) = T(y_1) = y_1$ or $T^2(x) = T(y_1) = x$, both of which leads to $T^n x = \{x,y_1\}\neq H$. As the same argument goes on, we have $k_1=1$, $k_2=2$, $\dots$, $k_m=m$. That is, we have proved that
\begin{align*}
y_1 = T(x),\quad y_2 = T(y_1), \quad \dots,\quad y_m = T(y_{m-1}).
\end{align*}

Now we prove by contradiction the last step: $x = T(y_m)$. Suppose $T(y_m)\neq x$. Then there exists $i\in \{1,\dots, m\}$ such that $T(y_m) = y_i$. As a consequence, one can easily check that $\{y_i,y_{i+1},\dots, y_m\}\subset H$ is an invariant set. Moreover, since $T(\{x,y_1,\dots, y_{i-1}\})=\{y_1,\dots, y_i\}$, we have 
\begin{align*}
T(H) = T \big( \{x,y_1,\dots, y_{i-1}\} \cup \{y_i,\dots, y_m\} \big) = T(\{x,y_1,\dots, y_{i-1}\})\cup T(\{y_i,\dots,y_m\}) = \{y_1,\dots, y_m\}\subset H,
\end{align*}
which contradicts the fact that $H$ is an invariant set. Therefore, there does not exist any $i\in \{1,\dots,m\}$ such that $T(y_m)=y_i$. That is, we have proved that $T(y_m)=x$ and thereby constructed the periodic motion as follows:
\begin{align*}
y_1 = T(x),\quad y_2 = T(y_1), \quad \dots, \quad y_m = T(y_{m-1}), \quad x = T(y_m).
\end{align*}
This concludes the proof. \QED

\begin{lemma}[Closures of invariant sets]\label{lem:inv-set-closure}
For any continuous map $T:\mathbb{R}^m \to \mathbb{R}^m$, the following statements hold:
\begin{enumerate}[label={\arabic*)}]
\item The closure of any positively invariant set of $T$ is positively invariant;
\item The closure of any bounded invariant set of $T$ is invariant. 
\end{enumerate}
\end{lemma}

\textbf{Proof:} Suppose $H$ is a positively invariant set of map $T$. Due to the continuity of $T$, if $T(x)\notin \overline{H}$ for some $x\in \partial H$, then there exists $y\in \intr(H)$ such that $T(y)\notin \overline{H}$. This implies that $H$ is not positively invariant, which leads to a contradiction. Therefore, $T(\partial H)\subset H$ and thus $T(H)\subset H$. This concludes the proof for statement~1).

According to statement~1), the closure of a bounded invariant set is positively invariant. Suppose $H$ is bounded and invariant. We have $T(\overline{H})\supset \overline{H}$, which leads to the following result: For any $x\in H$, there exists $y\in H$ such that $T(y)=x$. Since $H$ is bounded, $\partial H$ is well-defined and bounded. For any $x^*\in \partial H$, there exists a sequence $\{x_n\}\subset H$ such that $\{x_n\}\to x^*$. For any $x_n\in H$, there exists $y_n\in H$ such that $T(y_n)=x_n$. Now we obtain a sequence $\{y_n\}\subset H \subset \overline{H}$. Since $\overline{H}$ is a compact set, there exists a sequence $\{n_k\}\to \infty$ such that the subsequence $\{y_{n_k}\}$ converges to some $y^*\in \overline{H}$. Moreover, since map $T$ is continuous,
\begin{align*}
T(y^*) = T\big( \lim_{k\to \infty} y_{n_k} \big) = \lim_{k\to \infty} T(y_{n_k}) = \lim_{k\to\infty} x_{n_k} = x^* \in \overline{H}.
\end{align*}
Therefore, for any $x^*\in \overline{H}$, there exists $y^*\in \overline{H}$ such that $T(y^*)=x^*\in T(\overline{H})$. Now we have obtained both $T(\overline{H})\subset \overline{H}$ and $T(\overline{H})\supset \overline{H}$. Therefore, $T(\overline{H})=\overline{H}$, which concludes the proof of statement~2). \QED

\begin{lemma}[Invariant set and extension of motion]\label{lem:inv-set-extension-motion}
For any continuous map $T:\mathbb{R}^m\to \mathbb{R}^m$, a set $H$ is an invariant set of $T$ if and only if each motion $T^n x$ starting in $H$ has an extension in $H$. 
\end{lemma}

\textbf{Proof:} Suppose $H$ is an invariant set. Since $H$ is positively invariant, for any $x\in H$, we have $T^n(x)\in H$ for any $n\in \mathbb{N}$, which leads to $T^n\subset H$. On the other hand, since $H$ is also negatively invariant, for any $x\in H$, there exists $x_{-1}\in H$ such that $T(x_{-1})=x$. Let $T_{-1}(x)=x_{-1}$. Since $x_{-1}\in H$, there exists $x_{-2}\in H$ such that $T(x_{-2})\in H$. Following this argument and let $T_n(x)=T^n(x)$ for any $n\in \mathbb{N}$, we construct an extension $T_n x=\{T_n(x)\}_{n\in \mathbb{Z}}$ of the motion $T^n x$ and the extension $T_n x$ is in $H$.

Now suppose that each motion $T^n x$ starting in $H$ has an extension $T_n x$ in $H$. Since $x\in H$ leads to $T_n x \subset H$ for any $x\in H$, we have $T(x)\in H$, which implies $T(H)\subset H$. In addition, for any $x\in H$, since there exists $y=T_{-1}(x)\in H$ such that $T(y)=x$, we have $T(H)\supset H$. Therefore, $H$ is an invariant set of $T$. This concludes the proof.  \QED

\begin{lemma}[Properties of the largest invariant set]\label{lem:inv-set-largest}
For any continuous map $T\in \mathbb{R}^m \to \mathbb{R}^m$ and any set $E\subset \mathbb{R}^m$, if $M$ is the largest invariant set in $E$, then the following statements hold:
\begin{enumerate}[label=\arabic*)]
\item $M$ is the union all the extensions of motions that remain in $E$ for all $n\in \mathbb{Z}$;
\item $x\in M$ if and only if there exists an extension of motion $T_n x$ such that $T_n x\subset E$;
\item If $E$ is compact, then $M$ is compact.
\end{enumerate}
\end{lemma}

\textbf{Proof:} For any extension of motion staring at $x$, denoted by $T_n x$, one can easily check that the set $T_n x=\{T_n(x)\}_{n\in \mathbb{Z}}$ is invariant, i.e., $T(T_n x)=T_n x$. If $T_n x\in E$, then, according to Remark~\ref{remark:inv-set-connected-largest}, $T(M\cup T_n x)=T(M)\cup T(T_n x)=M\cup T_n x$ implies that $T_n x \subset M$. On the other hand, for any $x\in M$, since $M$ is invariant, we have $T(x)\in M$ and there exists $x_{-1}\in M$ such that $T(x_{-1})=x$. Then we further have $T^2(x)\in M$ and there exists $x_{-2}\in M$ such that $T(x_{-2})=x_{-1}$. Following this argument, we conclude that any $x\in M$ is in some extension of motion in $M$. This concludes the proof of statement~1).

Statement~2) is a straightforward result of statement~1).

Now we proceed to prove statement~3). For any sequence $\{x_n\}\subset M$ such that $\{x_n\}\to x^*$, since $E$ is compact, we have $x^*\in E$. Moreover, since $x_n\in M$ leads to $T^k(x_n)\in M$ for any $k\in \mathbb{N}$, due to the continuity of $T$, we have 
\begin{equation*}
\lim_{n\to \infty}T^k(x_n)=T^k\Big(\lim_{n\to \infty}x_n\Big)=T^k(x^*).
\end{equation*}
Therefore, the motion $T^n x^*$ satisfies $T^n x^* \subset E$. Moreover, since $M$ is invariant and $\{x_n\}\subset M$, for any $n\in \mathbb{N}$, there exists $y_n\in M$ such that $T(y_n)=x_n$. Now we obtain a sequence $\{y_n\}\subset M\subset E$. Since $E$ is compact, there exists a subsequence $\{y_{n_k}\}$ such that $\{n_k\}\to \infty$ and $\{y_{n_k}\}\to y^*\in E$ as $k\to \infty$. Due to the continuity of map $T$, we have
\begin{equation*}
T(y^*) = T\Big( \lim_{k\to \infty} y_{n_k} \Big) = \lim_{k\to \infty} T(y_{n_k})=\lim_{k\to \infty} x_{n_k} = x^*.
\end{equation*}
Let $T_{-1}(x^*)=y^*\in E$. Applying the same argument above for $x^*$ to $y^*$, we obtain $T_{-2}(x^*)\in E$ such that $T(T_{-2}(x^*))=T_{-1}(x^*)$. Continue this argument, we get an extension $T_n x^*$ of the motion $T^n x^*$ such that $T_n x^*\subset E$.According to statement~1), we have $x^*\in M$. Now we have proved that any $\{x_n\}\subset M$ such that $\{x_n\}\to x^*$ leads to $x^*\in M$. Therefore, $M$ is a compact set. \QED

\subsection{Properties of limit sets}

\begin{lemma}[Closed forms of limit sets]\label{lem:limit-set-closed-form}
Given a continuous map $T:\mathbb{R}^m \to \mathbb{R}^m$, for any $x\in \mathbb{R}^m$, the limit set of $x$ given by Definition~\ref{def:lim-set} satisfies
\begin{equation}\label{eq:limit-set-point-closed-form}
\Omega(x) = \cap_{j=0}^{\infty}\,\overline{\cup_{n=j}^{\infty} \{T^n(x)\}};
\end{equation}
Similarly, for any set $H\subset \mathbb{R}^m$, the limit set of $H$, given by Definition~\ref{def:lim-set}, satisfies
\begin{equation}\label{eq:limit-set-set-closed-form}
\Omega(H) = \cap_{j=0}^{\infty} \overline{\cup_{n=j}^{\infty} T^n(H)}.
\end{equation}
\end{lemma}

\textbf{Proof:} We first prove equation~\eqref{eq:limit-set-point-closed-form}. Suppose $y\in \Omega(x)$. For any $j\in \mathbb{N}$, since 
\begin{enumerate}
\item there exists $\{n_k\}_{k\in \mathbb{N}}$ such that $n_k\to \infty$ and $T^{n_k}(x)\to y$ as $k\to \infty$;
\item $\rho(y,\cup_{n=j}^{\infty}\{T^n(x)\})\le \rho(y,T^{n_k}(y))$ for any $k$ such that $n_k>j$,
\end{enumerate}
by letting $k\to \infty$, we have $\rho(y,\cup_{n=j}^{\infty}\{T^n(x)\})=0$, which implies that $y\in \overline{\cup_{n=j}^{\infty}\{T^n(x)\}}$. Since the argument above holds for any $j\in \mathbb{N}$, we have $y\in \cap_{j=1}^{\infty}\overline{\cup_{n=j}^{\infty}\{T^n(x)\}}$. On the other hand, for any $y\in \cap_{j=1}^{\infty}\overline{\cup_{n=j}^{\infty}\{T^n(x)\}}$, since $y\in \overline{\cup_{n=j}^{\infty}\{T^n(x)\}}$ for any $j\in \mathbb{N}$, we have $\rho(y,\cup_{n=j}^{\infty}\{T^n(x)\})=0$, that is, $\inf\limits_{z\in \cup_{n=j}^{\infty}\{T^n(x)\}}\rho(y,z)=0$, for any $j\in \mathbb{N}$. Let $0<\epsilon<1$ and $j=1$. There must exist some $n_1\ge 1$ such that $\rho(y,T^{n_1}(x))<\epsilon$. Then let $j=n_1+1$. There must exist some $n_2\ge j>n_1$ such that $\rho(y,T^{n_2}(x))<\epsilon^2$. Following this argument, we construct a subsequence $\{n_k\}_{k\in \mathbb{Z}_+}$ such that $\{n_k\} \to \infty$ and $\rho(y,T^{n_k}(x))<\epsilon^k$ for any $k\in \mathbb{Z}_+$, which implies that $T^{n_k}\to y$ as $k\to \infty$. Therefore, according to Definition~\ref{def:lim-set}, $y\in \Omega(x)$. This concludes the proof for equation~\eqref{eq:limit-set-point-closed-form}.

Similarly, equation~\eqref{eq:limit-set-set-closed-form} is proved by the following argument:
\begin{align*}
y\in \Omega(H) & \iff \exists\, \{n_k\} \to \infty \text{ and }\{y_k\}\subset H,\text{ s.t. }T^{n_k}(y_k)\to y\text{ as }k\to \infty, \\
               & \iff \forall\, j\in \mathbb{N},\text{ }\forall\, k\in \mathbb{N},\text{ }\exists\, n_k\ge j\text{ and }y_k\in H\text{ s.t. }T^{n_k}(y_k)\in \mathcal{B}(y,\epsilon^k), \text{ for some }0<\epsilon<1,\\ 
               & \iff \forall\, j\in \mathbb{N}, \text{ }\inf_{n\ge j,\,x\in H} \lVert y-T^n(x) \rVert=0, \\
               & \iff \forall\, j\in \mathbb{N}, \text{ }\rho( y,\cup_{n=j}^{\infty} T^n(H) )=0.
\end{align*}\QED

\begin{lemma}[Invariance and asymptotic properties of limit set $\Omega(x)$]\label{lem:limit-set-inv-asym}
For any continuous map $T:\mathbb{R}^{m}\to \mathbb{R}^m$ and any $x\in \mathbb{R}^m$, the following statements hold:
\begin{enumerate}[label=\arabic*)]
\item The limit set $\Omega(x)$ is closed and positively invariant;
\item If the motion $T^n x$ is a bounded set, then $\Omega(x)$ is
   \begin{enumerate}[label=\alph*)]
   \item nonempty;
   \item compact;
   \item invariant;
   \item invariantly connected;
   \item the smallest set that $T^n(x)$ approaches as $n\to \infty$.
   \end{enumerate}    
\end{enumerate}
\end{lemma}

\textbf{Proof:} The sketch of this proof can be found on Page 4 of~\cite{JPL:76}. Let $A_j(x)=\overline{\cup_{n=j}^{\infty}\{T^n(x)\}}$. According to Lemma~\ref{lem:limit-set-closed-form}, we have $\Omega(x)=\cap_{j=0}^{\infty}A_j(x)$. Each $A_j(x)$ is a closed set since it is the closure of a set of countably many points in $\mathbb{R}^m$. Therefore, as the intersection of countably many closed sets $A_j(x)$'s, $\Omega(x)$ is closed. Moreover, for any $y\in \Omega(x)$, by definition, there exists a sequence $\{n_k\}\to \infty$ such that $T^{n_k}(x)\to y$ as $k\to \infty$. Let $\ell_k=n_k+1$. Due to the continuity of map $T$, we have $T^{\ell_k}(x)=T\big(T^{n_k}(x)\big)\to T(y)$, which implies that $T(y)$ is also a limit point of $x$, i.e., $T(y)\in \Omega(x)$. Therefore, we obtain $T\big(\Omega(x)\big)\subset \Omega(x)$. This concludes the proof of statement~1). 

Now we prove statement~2). Since $T^n x$ is bounded, there exists a sequence $\{n_k\}\to \infty$ such that $T^{n_k}(x)$ converges. By definition, $\lim_{k\to \infty}T^{n_k}(x)\in \Omega(x)$. Therefore, $\Omega(x)$ is non-empty. This concludes the proof of statement~2)a).

By definition, $\Omega(x)\subset \overline{\cup_{n=0}^{\infty}\{T^n (x)\}}$. If $T^n x$ is bounded, then $\overline{\cup_{n=0}^{\infty}\{T^n (x)\}}$ is bounded, which in turn implied that $\Omega(x)$ is bounded. According to statement~1), $\Omega9x)$ is closed. Therefore, $\Omega(x)\subset \mathbb{R}^m$ is compact. This concludes the proof of statement~2)b).

For any $y\in \Omega(x)$, there exists $\{n_k\}\to \infty$ such that $T^{n_k}(x)\to y$. Let $\ell_k = n_k -1$, then $T\big( T^{\ell_k}(x) \big)\to y$ as $k\to \infty$. Since $\{T^{\ell_k}(x)\}_{k\in \mathbb{N}}\subset T^n x$ is bounded, there exists a sequence $\{k_r\}\to \infty$ such that $T^{\ell_{k_r}}(x)$ converges to some $z\in \Omega(x)$ as $r\to \infty$. Due to the continuity of $T$, we have
\begin{equation*}
T(z) = T\Big( \lim_{r\to \infty} T^{\ell_{k_r}}(x) \Big) = \lim_{r\to \infty} T^{\ell_{k_r}+1}(x)=\lim_{r\to \infty}T^{n_{k_r}}(x)=y.
\end{equation*}
Now we have shown that, for any $y\in \Omega(x)$, there exists $z\in \Omega(x)$ such that $T(z)=y$. Therefore, $T\big( \Omega(x) \big)\supset \Omega(x)$. We already have $T\big( \Omega(x) \big)\subset \Omega(x)$ according to statement~1). Therefore, $\Omega(x)$ is invariant. This concludes the proof of statement~2)c).

We prove statement~2)d) by contradiction. Suppose $\Omega(x)$ is a union of two disjoint sets $\Omega_1(x)$ and $\Omega_2(x)$, which are both closed, non-empty, and invariant. As a consequence, $d = \rho\big(\Omega_1(x),\Omega_2(x)\big)>0$. Due to the continuity of $T$ and $T\big( \Omega_1(x) \big)=\Omega_1(x)$, for any $0<\delta_1<d/2$, there exists $\epsilon_1$ satisfying $0<\epsilon_1<\delta_1$ such that $T\big( \mathcal{B}( \Omega_1(x),\epsilon_1 ) \big)\subset \mathcal{B}(\Omega_1(x),\delta_1)$. Here we adopt the generalized definition of an open ball around a set. That is, for any set $S$, $\mathcal{B}(S,\epsilon)=\{x\in \mathbb{R}^m\,|\,\rho(x,S)<\epsilon\}$. Similar to the argument above, we have that, for any $0<\delta_2<d/2$, there exists $0<\epsilon_2<\delta_2$ such that $T\big( \mathcal{B}(\Omega_2(x),\epsilon_2) \big)\subset \mathcal{B}(\Omega_2(x),\delta_2)$. Moreover, since $\Omega_1(x)\subset \Omega(x)$, there must exist a sequence $\{n_k\}\to \infty$ such that $T^{n_k}(x)$ enters $\mathcal{B}\big( \Omega_1(x),\epsilon_1 \big)$ for infinite times. For any $k$, there must exist a positive integer $\ell_k>n_k$ such that $T^{\ell_k}(x)\in \mathcal{B}\big( \Omega_2(x),\epsilon_2 \big)$, otherwise $\Omega_2(x)$ cannot be a subset of the limit set of $x$. However, since $T^{n_k}(x)\in \mathcal{B}\big( \Omega_1(x),\epsilon_1 \big)$, $T^{n_k+1}(x)$ cannot be in $\mathcal{B}\big( \Omega_2(x),\epsilon_2 \big)$. Therefore, for each $k\in \mathbb{N}$, there exists $\tilde{\ell}_k$ satisfying $n_k<\tilde{\ell}_k<\ell_k$ such that $T^{\tilde{\ell}_k}(x)\notin \mathcal{B}\big( \Omega_1(x),\epsilon_1 \big)$ and  $T^{\tilde{\ell}_k}(x)\notin \mathcal{B}\big( \Omega_2(x),\epsilon_2 \big)$. That is, for any $k\in \mathbb{N}$,
\begin{equation*}
T^{\tilde{\ell}_k}(x) \in \overline{T^n x}\,\setminus\, \Big( \mathcal{B}\big( \Omega_1(x),\epsilon_1 \big)\, \cup\, \mathcal{B}\big( \Omega_2(x),\epsilon_2 \big) \Big), 
\end{equation*} 
which is a bounded and closed set. Since there are infinitely many $\tilde{\ell}_k$'s such that $T^{\tilde{\ell}_k}(x)$ is in the compact set $\overline{T^n x}\,\setminus\, \Big( \mathcal{B}\big( \Omega_1(x),\epsilon_1 \big)\, \cup\, \mathcal{B}\big( \Omega_2(x),\epsilon_2 \big) \Big)$, this compact set must contain at least one limit point $z\in \Omega(x)$, which contradicts the assumption that $\Omega(x)=\Omega_1(x)\cup \Omega_2(x)$ since $z\notin \Omega_1(x)$ and $z\notin \Omega_2(x)$. This concludes the proof of statement~2)d).

To prove statement~2)e), we first prove by contradiction that $T^n(x)$ approaches $\Omega(x)$. We point out that $T^n (x)$ approaches $\Omega(x)$ if and only if $\rho\big( T^n(x),\Omega(x) \big)\to 0$ as $n\to \infty$. Suppose $\rho\big( T^n(x),\Omega(x) \big)$ does not converges to $0$. Then there exists $\epsilon>0$ such that, for any $k\in \mathbb{N}$, there exists $n\ge k$ such that $\rho\big( T^n(x),\Omega(x) \big)\ge \epsilon$. By letting $k=1,2,\dots$, we obtain a subsequence $\{n_k\}\to \infty$ such that $\rho\big( T^{n_k}(x),\Omega(x) \big)\ge \epsilon$ for any $k$. However, since $\{T^{n_k}(x)\}_{k\in \mathbb{N}}$ is bounded, there exists a sequence $\{k_r\}\to \infty$ such that $T^{n_{k_r}}\to y\in \Omega(x)$, which implies 
\begin{equation*}
\rho\big( T^{n_{k_r}}(x),y \big) \ge \rho\big( T^{n_{k_r}}(x),\Omega(x) \big) \to 0
\end{equation*} 
and contradicts with $\rho\big( T^{n_{k_r}}(x),\Omega(x) \big)\ge \epsilon$ for any $r$. Therefore $\rho\big( T^n(x),\Omega(x) \big)\to 0$.

Now we proceed to prove that $\Omega(x)$ is the smallest set that $T^n(x)$ approaches. For any set $A$, if $T^n(x)$ approaches $A$, then $T^n(x)$ also approaches $\overline{A}$. Therefore, we only need to discuss the case when $A$ is a closed set. Suppose $T^n(x)\to A$ as $n\to \infty$. Then we have that, for any sequence $\{n_k\}\to \infty$, $\rho\big( T^{n_k}(x),A \big)\to 0$. In addition, for any $y\in \Omega(x)$, there exists a sequence $\{n_k\}\to \infty$ such that $T^{n_k}\to y$ as $k\to \infty$. Suppose $y\notin A$. Since $A$ is closed, we have $\rho(A,y)>0$. For such sequence $\{n_k\}$ that $T^{n_k}(x)\to y$, we already obtain $\rho\big( T^{n_k}(x),A \big)\to 0$ and $\rho\big( T^{n_k}(x),y \big)\to 0$. Since $\rho(A,y)\le \rho\big( T^{n_k}(x),A \big) + \rho\big( T^{n_k}(x),y \big)$ for any $k\in \mathbb{N}$, we have $\rho(A,y)=0$, which contradicts $\rho(A,y)>0$. Therefore, $y\in A$ and thus $\Omega(x)\subset A$. This concludes the proof of statement~2)e).\QED 

The following lemma presents some important properties of the limit set of any set $H$, i.e., $\Omega(H)$. The proof follows the same line of argument in the proof of Lemma~\ref{lem:limit-set-inv-asym}.
\begin{lemma}[Invariance and asymptotic properties of limit set $\Omega(H)$]\label{lem:limit-set-H-inv-asym}
For any continuous map $T:\mathbb{R}^m \to \mathbb{R}^m$ and any set $H\subset \mathbb{R}^m$, the following statements for the limit set $\Omega(H)$ hold:
\begin{enumerate}[label=\arabic*)]
\item For any $x\in H$, $\Omega(x)\subset \Omega(H)$;
\item $\Omega(H)$ is closed and positively invariant;
\item If $\cup_{n=0}^{\infty} T^n(H)$ is bounded, then
   \smallskip
   \begin{enumerate}[label=\alph*)]
   \item $\Omega(H)$ is non-empty, compact, and invariant;
   \item for any $x\in H$, $T^n(x)$ approaches $\Omega(H)$;
   \item $\Omega(H)$ is the smallest set that $T^n(H)$ approaches as $n\to \infty$.   
   \end{enumerate}
\end{enumerate}
\end{lemma}

\begin{remark}
Unlike the limit set of a point, $\Omega(H)$ is not necessarily invariantly connected, even if $\cup_{n=0}^{\infty} T^n(H)$ is bounded. For example, consider the following one-dimension first-order difference equation system:
\begin{equation*}
x(n+1) = T\big( x(n) \big) = \epsilon x(n)^2 + (1-\epsilon)x(n),
\end{equation*} 
where $x(n)\in \mathbb{R}$ for any $n\in \mathbb{N}$ and $0<\epsilon<1$. One can easily check that the set of fixed points is $\{0,1\}$ and $\Omega([0,1])=\{0,1\}$, which is a union of two disjoint, non-empty, and closed invariant sets: $\{0\}$ and $\{1\}$.
\end{remark}

\begin{lemma}[Limit set of compact \& positively invariant set]\label{lem:limit-set-compact-pos-inv-K}
For any continuous map $T:\mathbb{R}^m\to \mathbb{R}^m$, suppose $E\subset \mathbb{R}^m$ is a compact and positively invariant set. The following statements hold:
\begin{enumerate}[label=\arabic*)]
\item $\Omega(E)=\cap_{n=0}^{\infty} T^n(E)$;
\item $\Omega(E)$ is non-empty, compact, and invariant;
\item $\Omega(E)$ is the largest invariant set in $E$. 
\end{enumerate}  
\end{lemma}

\textbf{Proof:} Since $E$ is compact and $T$ is continuous, $T^n(E)$ is compact for any $n\in \mathbb{N}$. Since $E$ is positively invariant, $T(E)\subset E$ and $T^{n+1}(E)\subset T^n(E)$ for any $n\in \mathbb{N}$. Therefore, by definition, 
\begin{equation*}
\Omega(E) = \cap_{j=0}^{\infty} \overline{\cup_{n=j}^{\infty}T^n(E)} = \cap_{j=0}^{\infty} \overline{T^j(E)} = \cap_{j=0}^{\infty} T^j(E).
\end{equation*}
This concludes the proof for statement~1).

Since $\cup_{n=0}^{\infty}T^n(E)=T^0(E)=E$ is bounded, according to statement~3)a) in Lemma~\ref{lem:limit-set-H-inv-asym}, $\Omega(E)$ is non-empty, compact, and invariant. This proves statement~2).

Now we prove by contradiction that $\Omega(E)$ is the largest invariant set in $E$. Suppose there exists a non-empty set $A\subset E$ such that
\begin{equation*}
A\cap \Omega(E)=\phi, \quad \text{and}\quad T\big( A\cup \Omega(E) \big)=A\cup \Omega(E).
\end{equation*}
Since $T\big( A\cup \Omega(E) \big)=T(A) \cup T\big( \Omega(E) \big)=T(A)\cup \Omega(E)$, we have $T(E)\supset T(A) \supset A$. In addition, 
\begin{equation*}
T\big( A\cup \Omega(E) \big) = A \cup \Omega(E) \quad\Rightarrow\quad T^n\big( A\cup \Omega(E) \big) = A\cup \Omega(E),\text{ for any }n\in \mathbb{N}. \quad\Rightarrow\quad \cap_{n=0}^{\infty} T^n\big( A\cup \Omega(E) \big)=A\cup \Omega(E).
\end{equation*} 
In the meanwhile,
\begin{align*}
T\big( A\cup \Omega(E) \big) = T(A) \cup \Omega(E) \quad & \Rightarrow \quad T^n\big( A\cup \Omega(E) \big)=T^n(A) \cup \Omega(E) \subset T^n(E) \cup \Omega(E)\\
& \Rightarrow \quad \cap_{n=0}^{\infty} T^n\big( A\cup \Omega(E) \big)\, \subset\, \cap_{n=0}^{\infty} \Big( T^n(E)\cup T^n\big( \Omega(E) \big) \Big) \,=\, \Big( \cap_{n=0}^{\infty} T^n(E) \Big) \cup \Omega(E)\,=\,\Omega(E)\\
& \Rightarrow\quad A\cup \Omega(E) = \Omega(E) \\
& \Rightarrow \quad A\subset \Omega(E),
\end{align*}  
which contradicts $A\cap \Omega(E)=\phi$. This concludes the proof for statement~3).\QED

\section{LaSalle Invariance Principle and Its Extension}
With all the preparation work in Section~1 and Section~2, now we are ready to present the main theorem on the original LaSalle invariance principle for discrete-time dynamical systems. The proof can be found on Page 6 of~\cite{JPL:76}.
\begin{theorem}[Discrete-time LaSalle invariance principle]\label{thm:classic-LaSalle}
Let $G$ be any set in $\mathbb{R}^m$. Consider a difference equations system defined by a map $T:\mathbb{R}^m \to \mathbb{R}^m$ that is well-defined for any $x\in \overline{G}$ and continuous at any $x\in G$. Suppose there exists a scalar map $V:\mathbb{R}^m \to \mathbb{R}$ satisfying
\begin{enumerate}[label=\roman*)]
\item $V(x)$ is continuous at any $x\in \overline{G}$;
\item $V\big( T(x) \big)-V(x)\le 0$ for any $x\in G$.
\end{enumerate} 
For any $x_0\in G$, if the solution to the following initial-value problem
\begin{equation*}
\begin{cases}
\displaystyle x(n+1) & = T\big( x(n) \big),\\
\displaystyle x(0) & = x_0
\end{cases}
\end{equation*}
satisfies that $\cup_{n=0}^{\infty}\{x(n)\}$ is bounded and $x(n)\in G$ for any $n\in \mathbb{N}$, then there exists $c\in \mathbb{R}$ such that $x(n)\to M\cap V^{-1}(c)$ as $n\to \infty$, where $V^{-1}(c)=\{x\in \mathbb{R}^m\,|\,V(x)=c\}$ and $M$ is the largest invariant set in $E=\big{\{} x\in \overline{G}\,\big|\, V\big( T(x) \big)-V(x)=0 \big{\}}$.
\end{theorem}

\textbf{Proof:} Let $X=\cup_{n=0}^{\infty}\{x(n)\}$. We have that $\overline{X}\subset \overline{G}$ is compact. Since $V(x)$ is continuous on $\overline{G}$, $V(x)$ is lower bounded on $\overline{X}$. Moreover, since $V\big(x(n)\big)$ is non-increasing with $n$ for any $n\in \mathbb{N}$, there exists $c\in \mathbb{R}$ such that $V\big( x(n) \big)\to c$ as $n\to \infty$.

For any $y\in \Omega(x_0)$, since there exists $\{n_k\}\to \infty$ such that $x(n_k)=T^{n_k}(x_0)\to y$ as $k\to \infty$ and due to the continuity of $V$, we have $V\big( x_{n_k} \big)\to V(y)$. Therefore,
\begin{equation*}
V(y) = \lim_{k\to \infty} V\big( x(n_k) \big) = \lim_{n\to \infty} V\big( x(n) \big) = c, \quad \text{ for any }y\in \Omega(x_0),
\end{equation*}
which leads to $\Omega(x_0)\in V^{-1}(c)$.

Moreover, for any $y\in \Omega(x_0)$, since $X$ is bounded, according to Lemma~\ref{lem:limit-set-inv-asym}, $\Omega(x_0)$ is invariant. Therefore, $T(y)\in \Omega(x_0)$, which implies that $V\big( T(y) \big)=c$ and thus $V\big( T(y) \big)-V(y)=0$ for any $y\in \Omega(x_0)$. Now we obtain $\Omega(x_0)\subset E$. Since $\Omega(x_0)\subset E$ and $M$ is the largest invariant set of $E$, we have $\Omega(x_0)\subset M\cap V^{-1}(c)$. Finally, since $T^n(x)$ approaches $\Omega(x_0)$ as $n\to \infty$, $x(n)\to M\cap V^{-1}(c)$ as $n\to \infty$. \QED

The classic LaSalle invariance principle stated in Theorem~\ref{thm:classic-LaSalle} requires that both $T$ and $V$ are well-defined on $\overline{G}$. Below we present a simple extension of the classic LaSalle invariance principle. This extension establish the converge of the solution $x(n)$ when $T$ and $V$ are not defined on $\partial G$ but $x(n)$ is uniformly bounded from $\partial G$ after some finite time $N$.

\begin{theorem}[Extension of LaSalle invariance principle]\label{thm:LaSalle-extension}
Consider the following difference equations system:
\begin{equation}\label{eq:thm-LaSalle-extension-system}
x(n+1) = T\big( x(n) \big),
\end{equation}
where $T:\mathbb{R}^m \to \mathbb{R}^m$ is continuous on some set $G\subset \mathbb{R}^m$. Suppose there exists a map $V:\mathbb{R}^m \to \mathbb{R}$ satisfying: i) $V(x)$ is continuous at any $x\in G$; ii) $V\big( T(x) \big)-V(x)\le 0$ for any $x\in G$. For any $x_0\in G$, if there exists a compact set $G_c$ and $N\in \mathbb{N}$ such that the solution $x(n)$ to equation~\eqref{eq:thm-LaSalle-extension-system} with $x(0)=x_0$ satisfies $x(n)\in G_c$ for any $n\ge N$, then there exists $c\in \mathbb{R}$ such that $x(n)\to M\cap V^{-1}(c)$ as $n\to \infty$, where $V^{-1}(c) = \{x\in \mathbb{R}^m\,|\,V(x)=c\}$ and $M$ is the largest invariant set in $E=\big{\{}x\in G_c\,|\,V\big( T(x) \big)-V(x)=0\big{\}}$.  
\end{theorem}

\textbf{Proof:} For any such $x_0\in G$, since $x(n)$ is the solution, we have $x(n)=T^n(x_0)$ for any $n\in \mathbb{N}$. According to Lemma~\ref{lem:limit-set-closed-form}, 
\begin{equation*}
\Omega(x_0) \,=\, \cap_{j=0}^{\infty} \overline{\cup_{n=j}^{\infty}\{T^n( x_0 )\}}\,\subset\,\overline{\cup_{n=N}^{\infty}\{T^{n}( x_0 )\}}\,\subset \, G_c.
\end{equation*}
For any $n\ge N$, since $x(n)\in \overline{\cup_{j=N}^{\infty}\{x(j)\}}\subset G_c$ and $V$ is continuous on $G_c$, $V\big( x(n) \big)$ is uniformly lower bounded for all $n\ge N$. In addition, since $V\big( x(n) \big)$ is non-increasing, there exists $c\in \mathbb{R}$ such that $\lim_{n\to \infty}V\big( x(t) \big)=c$. 

For any $y\in \Omega(x_0)$, there exists a sequence $\{n_k\}\to \infty$ such that $\{x(n_k)\}\to y$ as $k\to \infty$. Since $V$ is continuous on $G$, we have $V\big( x(n_k) \big) \to V(y)$. Moreover, $\lim_{n\to \infty} V(n)=c$ leads to $V(y)=c$. Therefore, $\Omega(x_0)\subset V^{-1}(c)$.

According to Lemma~\ref{lem:limit-set-inv-asym}, $\Omega(x_0)$ is invariant. Therefore, $T(y)\in \Omega(x_0)$ for any $y\in \Omega(x_0)$, which in turn implies that $V\big( T(y) \big)-V(y)=0$ for any $y\in \Omega(x_0)$ and thereby $\Omega(x_0)\subset E$. Moreover, since $\Omega(x_0)$ is invariant, $\Omega(x_0)\subset M$. Therefore, $\Omega(x_0)\subset M\cap V^{-1}(c)$. Since $x(n)=T^n(x)$ approaches $\Omega(x_0)$, $x(n)\to M\cap V^{-1}(c)$ as $n\to \infty$. This concludes the proof.
\QED

\section{Advanced Versions of LaSalle Invariance Principle}
In this section we provide an incomplete list of references on the extensions and more advanced versions of LaSalle invariance principle for the interest of further reading. In Section~8, Chapter~1 of~\cite{JPL:76}, the author discusses the vector Lyapunov functions; Cort\'{e}s et al.~\cite{JC-SM-TK-FB:02j} propose a sufficient condition for the convergence of discrete-time systems to the fixed points. Hale~\cite{JKH:69-article} extends the LaSalle invariance principles to autonomous systems with infinite dimensions; Shevitz and Paden~\cite{DS-BP:94} discusses LaSalle invariance principle in non-smooth systems; Extensions of LaSalle's results to switched systems are provided by Hespanha et al.~\cite{JH-DL-DA-EDS:05}, Bacciotti et al.~\cite{AB-LM:05}, and Mancilla et al.~\cite{JLMA-RAG:06}; Alberto et al.~\cite{LFCA-TRC-ACPM:07} consider the invariance principle for discrete-time dynamical systems with generalized Lyapunov functions of which the first difference are positive in some bounded regions; Results on Lyapunov functions and invariance principles for difference inclusions systems can be found in the research articles by Kellett and Teel~\cite{CMK-ART:04} and Bullo et al.~\cite{FB-RC-PF:08u} (see Lemma 4.1), as well as in the book by Bullo et al.~\cite{FB-JC-SM:09} (see Theorem 1.21); We refer to the book by Goebel et al.~\cite{RG-RGS-ART:12} for a systematic treatment of hybrid dynamical systems.

\bibliographystyle{plainurl}
\bibliography{alias,FB,Main}

\begin{thebibliography}{10}

\bibitem{LFCA-TRC-ACPM:07}
L.~F.~C. Alberto, T.~R. Calliero, and A.~C.~P. Martins.
\newblock An invariance principle for nonlinear discrete autonomous dynamical
  systems.
\newblock {\em IEEE Transactions on Automatic Control}, 52:692--697, 2007.
\newblock \href {http://dx.doi.org/10.1109/TAC.2007.894532}
  {\path{doi:10.1109/TAC.2007.894532}}.

\bibitem{AB-LM:05}
A.~Bacciotti and L.~Mazzi.
\newblock An invariance principle for nonlinear switched systems.
\newblock {\em Systems \& Control Letters}, 54(11):1109--1119, 2005.
\newblock \href {http://dx.doi.org/10.1016/j.sysconle.2005.04.003}
  {\path{doi:10.1016/j.sysconle.2005.04.003}}.

\bibitem{FB-RC-PF:08u}
F.~Bullo, R.~Carli, and P.~Frasca.
\newblock Gossip coverage control for robotic networks: {D}ynamical systems on
  the space of partitions.
\newblock {\em SIAM Journal on Control and Optimization}, 50(1):419--447, 2012.
\newblock \href {http://dx.doi.org/10.1137/100806370}
  {\path{doi:10.1137/100806370}}.

\bibitem{FB-JC-SM:09}
F.~Bullo, J.~Cort{\'e}s, and S.~Mart{\'\i}nez.
\newblock {\em Distributed Control of Robotic Networks}.
\newblock Princeton University Press, 2009.
\newblock URL: \url{http://www.coordinationbook.info}.

\bibitem{JC-SM-TK-FB:02j}
J.~Cort{\'e}s, S.~Mart{\'\i}nez, T.~Karatas, and F.~Bullo.
\newblock Coverage control for mobile sensing networks.
\newblock {\em IEEE Transactions on Robotics and Automation}, 20(2):243--255,
  2004.
\newblock \href {http://dx.doi.org/10.1109/TRA.2004.824698}
  {\path{doi:10.1109/TRA.2004.824698}}.

\bibitem{RG-RGS-ART:12}
R.~Goebel, R.~G. Sanfelice, and A.~R. Teel.
\newblock {\em Hybrid Dynamical Systems: Modeling, Stability, and Robustness}.
\newblock Princeton University Press, 2012.

\bibitem{JKH:69-article}
J.~K. Hale.
\newblock Dynamical systems and stability.
\newblock {\em Journal of Mathematical Analysis and Applications}, 26:39--59,
  1969.
\newblock \href {http://dx.doi.org/10.1016/0022-247X(69)90175-9}
  {\path{doi:10.1016/0022-247X(69)90175-9}}.

\bibitem{JH-DL-DA-EDS:05}
J.~Hespanha, D.~Liberzon, D.~Angeli, and E.~D. Sontag.
\newblock Nonlinear norm-observability notions and stability of switched
  systems.
\newblock {\em IEEE Transactions on Automatic Control}, 50(2):154--168, 2005.
\newblock \href {http://dx.doi.org/10.1109/TAC.2004.841937}
  {\path{doi:10.1109/TAC.2004.841937}}.

\bibitem{CMK-ART:04}
C.~M. Kellett and A.~R. Teel.
\newblock Smooth {Lyapunov} functions and robustness of stability for
  difference inclusions.
\newblock {\em Systems \& Control Letters}, 52:395--405, 2004.
\newblock \href {http://dx.doi.org/10.1016/j.sysconle.2004.02.015}
  {\path{doi:10.1016/j.sysconle.2004.02.015}}.

\bibitem{HKK:02}
H.~K. Khalil.
\newblock {\em Nonlinear Systems}.
\newblock Prentice Hall, 3 edition, 2002.

\bibitem{JPL:76}
J.~P. LaSalle.
\newblock {\em The Stability of Dynamical Systems}.
\newblock SIAM, 1976.
\newblock \href {http://dx.doi.org/10.1137/1.9781611970432}
  {\path{doi:10.1137/1.9781611970432}}.

\bibitem{JLMA-RAG:06}
J.~L. Mancilla-Aguilar and R.~A. Garc\'{i}a.
\newblock An extension of {LaSalle's} invariance principle for switched
  systems.
\newblock {\em Systems \& Control Letters}, 55:376--384, 2006.
\newblock \href {http://dx.doi.org/10.1016/j.sysconle.2005.07.009}
  {\path{doi:10.1016/j.sysconle.2005.07.009}}.

\bibitem{DS-BP:94}
D.~Shevitz and B.~Paden.
\newblock Lyapunov stability theory of nonsmooth systems.
\newblock {\em IEEE Transactions on Automatic Control}, 39(9):1910--1914, 1994.

\bibitem{MV:02}
M.~Vidyasagar.
\newblock {\em Nonlinear Systems Analysis}.
\newblock SIAM, 2002.
\newblock \href {http://dx.doi.org/10.1137/1.9780898719185}
  {\path{doi:10.1137/1.9780898719185}}.

\end{thebibliography}
\end{document}